\documentclass[11pt]{amsart}
\usepackage{amsthm,amssymb,latexsym,graphicx,comment}
\usepackage{epstopdf}
\def\Re{\text{\rm Re}\,}

\def\dis{\displaystyle}

\makeatletter
\@namedef{subjclassname@2010}{%F
  \textup{2010} Mathematics Subject Classification}
\makeatother

\newtheorem{thm}{Theorem}
\newtheorem{lem}[thm]{Lemma}

\begin{document}

\baselineskip=15pt

\title[Suita Conjecture for Some Convex Ellipsoids in
$\mathbb C^2$]
{On the Suita Conjecture \\ for Some Convex Ellipsoids in
$\mathbb C^2$}

\thanks {The first named author was supported by the Ideas Plus 
grant 0001/ID3/2014/63 of the Polish Ministry Of Science and 
Higher Education and the second named author by 
the Polish National Science Centre grant  
2011/03/B/ST1/04758}

\def\nl{\newline\phantom{a}\hskip 7pt}

\author{Zbigniew B\l ocki, W\l odzimierz Zwonek}
\address{Uniwersytet Jagiello\'nski \nl Instytut Matematyki
\nl \L ojasiewicza 6 \nl 30-348 Krak\'ow \nl Poland
\nl {\rm Zbigniew.Blocki@im.uj.edu.pl 
\nl Wlodzimierz.Zwonek@im.uj.edu.pl}}

\begin{abstract}
It has been recently shown that for a convex domain $\Omega$
in $\mathbb C^n$ and $w\in\Omega$ the function 
$F_\Omega(w):=\big(K_\Omega(w)\lambda(I_\Omega(w))\big)^{1/n}$,
where $K_\Omega$ is the Bergman kernel on the diagonal and
$I_\Omega(w)$ the Kobayashi indicatrix, satisfies 
$1\leq F_\Omega\leq 4$. While the lower bound is optimal, not
much more is known about the upper bound. In general it is quite 
difficult to compute $F_\Omega$ even numerically and the highest 
value of it obtained so far is $1.010182\dots$ In this paper we present
precise, although rather complicated formulas for the ellipsoids
$\Omega=\{|z_1|^{2m}+|z_2|^2<1\}$ (with $m\geq 1/2$) and all $w$,  
as well as for $\Omega=\{|z_1|+|z_2|<1\}$ and $w$ on the
diagonal. The Bergman kernel for those ellipsoids had been
known, the main point is to compute the volume of the Kobayashi
indicatrix. It turns out that in the second case the function 
$\lambda(I_\Omega(w))$ is not $C^{3,1}$.
\end{abstract}

\makeatletter
\@namedef{subjclassname@2010}{%
  \textup{2010} Mathematics Subject Classification}
\makeatother

%\subjclass[2010]{}

%\keywords{}

\maketitle

\section*{\bf Introduction}

For a convex domain $\Omega$ in $\mathbb C^n$ and $w\in\Omega$ 
the following estimates have been recently established:
\begin{equation}\label{est}
\frac 1{\lambda(I_\Omega(w))}\leq K_\Omega(w)\leq
\frac {4^n}{\lambda(I_\Omega(w))}.
\end{equation}
Here
  $$K_\Omega(w)=\sup\{|f(w)|^2\colon f\in\mathcal O(\Omega),\
     \int_\Omega|f|^2d\lambda\leq 1\}$$
is the Bergman kernel on the diagonal and
  $$I_\Omega(w)=\{\varphi'(0)\colon \varphi\in\mathcal O(\Delta,\Omega),
     \ \varphi(0)=w\}$$
is the Kobayashi indicatrix, where $\Delta$ denotes the unit disc.
The first inequality in \eqref{est} was shown in \cite{B2}, the proof
uses $L^2$-estimates for $\bar\partial$ and Lempert's theory \cite{L}.
It is optimal, for example if $\Omega$ is balanced with respect to
$w$ (that is every intersection of $\Omega$ with a complex line 
containing $w$ is a disc) then we have equality. It can be viewed
as a multi-dimensional version of the Suita conjecture \cite{S}
proved in \cite{B1} (see also \cite{GZ} for the precise characterization
when equality holds).

The second equality in \eqref{est}
was proved in \cite{BZ} using rather elementary methods. It was
also shown that the constant $4$ can be replaced by 
$16/\pi^2=1.6211\dots$ if $\Omega$ is in addition symmetric with respect
to $w$. We can write \eqref{est} as
  $$1\leq F_\Omega(w)\leq 4,$$
where $F_\Omega(w):=\big(K_\Omega(w)\lambda(I_\Omega(w))\big)^{1/n}$
is a biholomorphically invariant function in $\Omega$. It is not clear 
what the optimal upper bound should be. It was in fact quite difficult
to prove that one can at all have $F_\Omega>1$. It was done in \cite{BZ}
for ellipsoids of the form $\{|z_1|+|z_2|^{2m}+\dots+\dots|z_n|^{2m}<1\}$,
where $m\geq 1/2$ and $w=(b,0,\dots,0)$. The function $F_\Omega$ was also 
computed numerically for the ellipsoid $\Omega=\{|z_1|^{2m}+|z_2|^2<1\}$, 
$m\geq 1/2$, based on an implicit formula for the Kobayashi function from 
\cite{6a}. Our first result is the precise formula in this case:

\begin{thm} \label{thm1}
For $m\geq 1/2$ define 
  $$\Omega_m=\{z\in\mathbb C^2\colon |z_1|^{2m}+|z_2|^2<1\}.$$
Then for $m\neq 2/3$, $m\neq 2$ and $b$ with $0\leq b<1$, we have
  $$\begin{aligned}
  \lambda(I_{\Omega_m}((b,0)))
      =\pi^2&\left[-\frac{m-1}{2m(3m-2)(3m-1)}b^{6m+2}
     -\frac{3(m-1)}{2m(m-2)(m+1)}b^{2m+2}\right.\\
    &\ \ \ \ \ \ \ \ \ \left.+\frac m{2(m-2)(3m-2)}b^6
    +\frac{3m}{3m-1}b^4-\frac{4m-1}{2m}b^2+\frac m{m+1}\right].
     \end{aligned}$$
For $m=2/3$ and $m=2$ one has 
  $$\lambda(I_{\Omega_{2/3}}((b,0)))=
    \frac{\pi^2}{80}\left(-65b^6+40b^6 \log b+160b^4
       -27b^{10/3}-100b^2+32\right),$$
  $$\lambda(I_{\Omega_2}((b,0)))
    =\frac{\pi^2}{240}\left(-3b^{14}-25b^6-120b^6\log b
       +288b^4-420b^2+160\right).$$
\end{thm}

The general formula for the Kobayashi function for $\Omega_m$  
is known, see \cite{6a}, but it is implicit in the sense that it requires 
solving a nonlinear equation which is polynomial of degree $2m$ if 
it is an integer. It turns out however that the volume of the Kobayashi 
indicatrix for $\Omega_m$, that is the set where the Kobayashi function 
is not bigger than 1, can be found explicitly. It would be interesting
to check whether Theorem \ref{thm1} also holds in the non-convex case,
that is when $0<m<1/2$ (see \cite{PZ} for computations of the
Kobayashi metric in this case). 

The formula for the Bergman kernel for this ellipsoid is well 
known (see e.g. \cite{JP}, Example 6.1.6):
  $$K_{\Omega_m}(w)=\frac 1{\pi^2}(1-|w_2|^2)^{1/m-2}
     \frac{(1/m+1)(1-|w_2|^2)^{1/m}+(1/m-1)|w_1|^2}
       {\big((1-|w_2|^2)^{1/m}-|w_1|^2\big)^3},$$
so that
  $$K_{\Omega_m}((b,0))=\frac{m+1+(1-m)b^2}{\pi^2m(1-b^2)^3},$$
and we can obtain the following graphs of $F_{\Omega_m}((b,0))$ for 
example for $m=4$, 8, 16, 32, 64 and 128:

\vskip 10pt

\centerline{\includegraphics[scale=1]{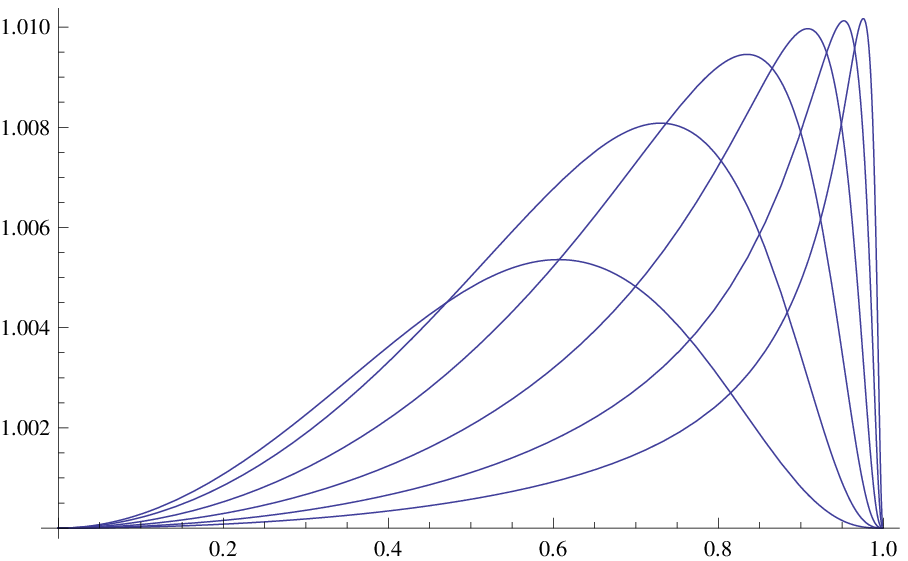}}

\vskip 10pt

They are consistent with the graphs from \cite{BZ} obtained
numerically using the implicit formula from \cite{6a}. Note that 
for $t\in\mathbb R$ and $a\in\Delta$ the mapping
  $$\Omega_m\ni z\longmapsto\left(e^{it}\frac{(1-|a|^2)^{1/2m}}
     {(1-\bar az_2)^{1/m}}z_1,\frac{z_2-a}{1-\bar az_2}\right)$$
is a holomorphic automorphism of $\Omega_m$ and therefore 
$F_{\Omega_m}((b,0))$ where $0\leq b<1$ attains all values of 
$F_{\Omega_m}$ in $\Omega_m$. One can show numerically that
  $$\sup_{m\geq 1/2}\sup_{\Omega_m}F_{\Omega_m}=1.010182\dots$$
which was already noticed in \cite{BZ}. This is the highest 
value of $F_\Omega$ (in arbitrary dimension) obtained so far.

In \cite{BZ} it was also shown that for $\Omega=\{|z_1|+|z_2|<1\}$ 
and $b$ with $0<b<1$ one has 
  $$\lambda(I_\Omega((b,0))=\frac{\pi^2}6(1-b)^4\big((1-b)^4+8b\big),$$
so that in particular similarly as in Theorem \ref{thm1}
it is an analytic function on this part of $\Omega$. 
This raises a question whether $\lambda(I_\Omega(w))$
is smooth in general. In \cite{BZ} it was also predicted that 
the highest value of $F_\Omega$ for convex $\Omega$ in 
$\mathbb C^2$ should be attained for 
for $\Omega=\{|z_1|+|z_2|<1\}$ on the diagonal. 
The following result will answer both of these questions
in the negative:

\begin{thm}\label{thm2}
Let $\Omega=\{z\in\mathbb C^2\colon |z_1|+|z_2|<1\}$. 
Then for $b$ with $0\leq b\leq 1/4$ we have  
\begin{equation}\label{b14}
  \lambda(I_\Omega((b,b)))=\frac{\pi^2}6
      \big(30b^8-64b^7+80b^6-80b^5+76b^4-16b^3-8b^2+1\big)
\end{equation}
and when $1/4\leq b<1/2$
\begin{equation}\label{a14}
\begin{aligned}
  \lambda(I_\Omega&((b,b)))=
    \frac{2\pi^2b(1-2b)^3\left(-2b^3+3b^2-6b+4\right)}{3 (1-b)^2}\\
    &+\frac{\pi\left(30 b^{10}-124 b^9+238 b^8-176 b^7-
     260 b^6+424 b^5-76 b^4-144 b^3+89b^2-18 b+1\right)}{6(1-b)^2}\\
    & \ \ \ \ \ \ \ \ \ \times\arccos\left(-1+\frac{4b-1}{2 b^2}\right)\\
    &+\frac{\pi(1-2 b)\left(-180 b^7+444 b^6-554 b^5+754 b^4-1214 b^3+922
       b^2-305 b+37\right)}{72 (1-b)}\sqrt{4 b-1}\\
    &+\frac{4\pi b(1-2b)^4\left(7 b^2+2 b-2\right)}{3(1-b)^2}
         \arctan\sqrt{4b-1}\\
    &+\frac{4\pi b^2(1-2 b)^4(2-b)}{(1-b)^2}
        \arctan\frac{1-3 b}{(1-b) \sqrt{4b-1}}.
\end{aligned}
\end{equation}
The function
  $$b\longmapsto\lambda(I_\Omega((b,b)))$$
is $C^3$ on the interval $(0,1/2)$ but not $C^{3,1}$ at $1/4$.
\end{thm}

Again, the formula for the Bergman metric for this ellipsoid 
is known, see \cite{HP} or \cite{JP}, Example 6.1.9:
  $$K_\Omega(w)=\frac 2{\pi^2}\cdot\frac{3(1-|w|^2)^2(1+|w|^2)
              +4|w_1|^2|w_2|^2(5-3|w|^2)}
     {\big((1-|w|^2)^2-4|w_1|^2|w_2|^2\big)^3},$$
so that
\begin{equation}\label{bergman}
  K_\Omega((b,b))=\frac{2(3-6b^2+8b^4)}{\pi^2(1-4b^2)^3}.
\end{equation}
The first part of Theorem \ref{thm2}, formula \eqref{b14}
on the interval $(0,1/4)$, is easier to prove than the
second one. Combining it with \eqref{bergman} one obtains 
the following graph of $F_\Omega((b,b))$ for $b\in(0,1/4)$:

\vskip 10pt

\centerline{\includegraphics[scale=1]{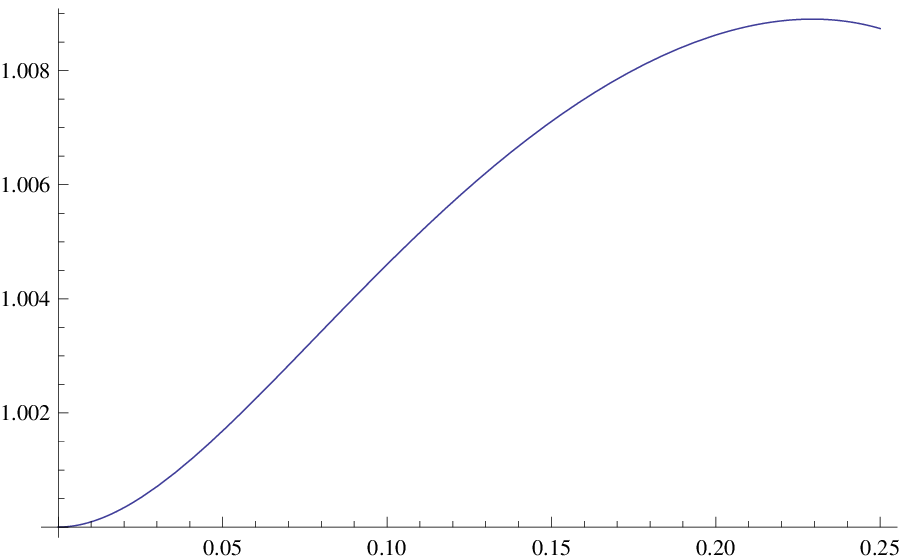}}

\vskip 10pt

\noindent One can show that its analytic continuation to 
$(0,1/2)$ attains values below 1 and thus it follows already
from \eqref{est} that $F_\Omega$ cannot be analytic. To conclude
that it is in fact not $C^{3,1}$ one has to prove much harder 
formula \eqref{a14}. Here is the full picture on the interval 
$(0,1/2)$, the analytic continuation of $F_\Omega$ from $(0,1/4)$ 
and the actual graph of $F_\Omega$:

\vskip 10pt

\centerline{\includegraphics[scale=1]{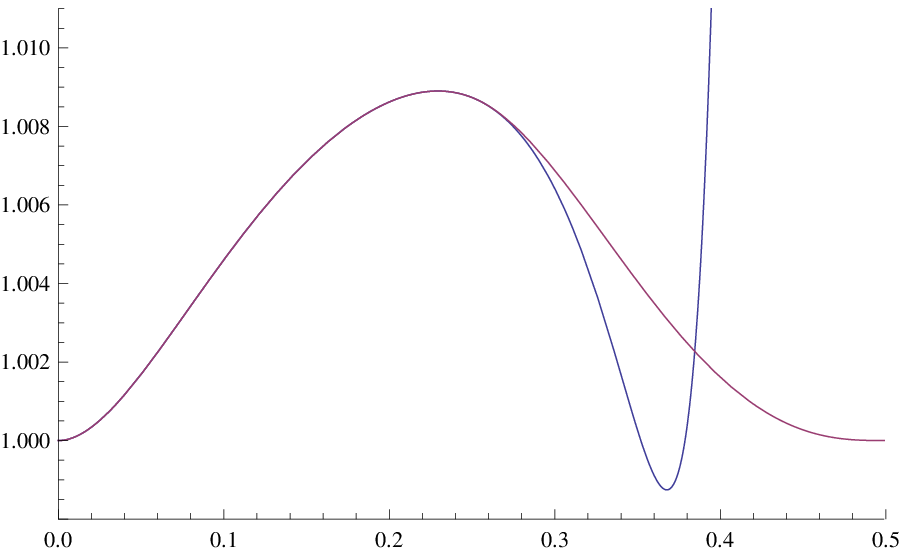}}

\vskip 10pt

\noindent One can check that the maximal value of $F_\Omega((b,b))$
for $b\in (0,1/2)$ is $1.008902\dots$

All pictures and numerical computations in this paper, as well as a lot 
of formal ones in the proofs of Theorems \ref{thm1} and \ref{thm2} have
been done using {\it Mathematica}.

\section{\bf General formula for geodesics in convex complex ellipsoids}

Boundary of the Kobayashi indicatrix of a convex domain $\Omega$ at $w$ 
consists of the vectors $\varphi'(0)$ where $\varphi\in\mathcal O(\Delta,\Omega)$
is a geodesic of $\Omega$ satisfying $\varphi(0)=w$. Theorems \ref{thm1} and 
\ref{thm2} will be proved using a general formula for geodesics in convex 
complex ellipsoids from \cite{JPZ} based on Lempert's theory \cite{L} 
describing geodesics of smooth strongly convex domains.

For $p=(p_1,\dots,p_n)$ with $p_j\geq 1/2$ set
  $$\mathcal E(p)
     =\{z\in\mathbb C^n: |z_1|^{2p_1}+\dots+|z_n|^{2p_n}<1\}$$
and $A\subset\{1,\dots,n\}$ define
\begin{equation*}\label{geod}
  \varphi_j(\zeta)
     =\begin{cases}\dis a_j\frac{\zeta-\alpha_j}{1-\bar\alpha_j\zeta}
       \left(\frac{1-\bar\alpha_j\zeta}
           {1-\bar\alpha_0\zeta}\right)^{1/p_j},\ &j\in A\\ \\
      \dis a_j\left(\frac{1-\bar\alpha_j\zeta}
         {1-\bar\alpha_0\zeta}\right)^{1/p_j},
              &j\notin A\end{cases},
\end{equation*}
where $a_j\in\mathbb C_\ast$, $\alpha_0,\alpha_j\in\Delta$ 
for $j\in A$, $\alpha_j\in\bar\Delta$ for $j\notin A$,
\begin{equation}\label{c1}
\alpha_0=|a_1|^{2p_1}\alpha_1+\dots+|a_n|^{2p_n}\alpha_n,
\end{equation}
and
\begin{equation}\label{c2}
1+|\alpha_0|^2=|a_1|^{2p_1}(1+|\alpha_1|^2)+\dots
     +|a_n|^{2p_n}(1+|\alpha_n|^2).
\end{equation}
A component $\varphi_j$ has a zero in $\Delta$ if and only if $j\in A$. 
We have
\begin{equation}\label{w0}
  \varphi_j(0)
     =\begin{cases}\dis -a_j\alpha_j,\ &j\in A\\
      \dis a_j, &j\notin A\end{cases},
\end{equation}
and
\begin{equation}\label{w1}
  \varphi_j'(0)
     =\begin{cases}\dis a_j\left(1+\big(\frac 1{p_j}-1\big)|
        \alpha_j|^2-\frac{\alpha_j\bar\alpha_0}{p_j}\right),\ 
            &j\in A\\ \\
      \dis a_j\frac{\bar\alpha_0-\bar\alpha_j}{p_j},
              &j\notin A\end{cases}.
\end{equation}
For $w\in\mathcal E(p)$ the set of vectors $\varphi'(0)$ where 
$\varphi(0)=w$ forms a subset of 
$\partial I^K_{\mathcal E(p)}(w)$ of a full measure. The geodesics
in $\mathcal E(p)$ are uniquely determined: for a given 
$w\in\mathcal E(p)$ and $X\in(\mathbb C^n)_\ast$ there exists 
unique geodesic $\varphi\in\mathcal O(\Delta,\mathcal E(p))$ 
such that $\varphi(0)=w$ and $\varphi'(0)=X$.

\section{\bf Proof of Theorem \ref{thm1}}

First note that the formulas for $m=2/3$ and $m=2$ easily follow
from the first one by approximation. 
For $\Omega_m=\mathcal E(m,1)$ and $w=(b,0)$ there are two possibilities
for a geodesic $\varphi$:
either $\varphi$ crosses the axis $\{z_1=0\}$ or it does not. By 
$I_{12}$ and $I_2$ denote the respective parts of $I_{\Omega_m}(w)$.
In the first case $\varphi$ must be of the form
  $$\varphi(\zeta)=\left(a_1\frac{\zeta-\alpha_1}{1-\bar\alpha_1\zeta}
         \left(\frac{1-\bar\alpha_1\zeta}{1-\bar\alpha_0\zeta}\right)^{1/m},
         a_2\frac{\zeta-\alpha_2}{1-\bar\alpha_0\zeta}\right),$$
where $a_1,a_2\in\mathbb C_\ast$ and $\alpha_0,\alpha_1,\alpha_2\in\Delta$
satisfy \eqref{c1}, \eqref{c2}.
By \eqref{w0} and since $\varphi(0)=(b,0)$ we have $a_1=-b/\alpha_1$, 
$\alpha_2=0$ and by \eqref{c1} $\alpha_0=b^{2m}\alpha_1/|\alpha_1|^{2m}$. 
By \eqref{c2}
  $$1+b^{4m}|\alpha_1|^{2-4m}=b^{2m}|\alpha_1|^{-2m}
       \left(1+|\alpha_1|^2\right)+|a_2|^2,$$
that is
\begin{equation}\label{a2}
  |a_2|^2=(1-b^{2m}|\alpha_1|^{-2m})(1-b^{2m}|\alpha_1|^{2-2m}).
\end{equation}
Since $\alpha_0,\alpha_1\in\Delta_\ast$, it follows that $b<|\alpha_1|<1$. 
Write $\alpha_1=-re^{-it}$, $a_2=\rho e^{is}$, then 
by \eqref{w1} and \eqref{a2}, with $b<r<1$,
  $$\begin{aligned}\varphi'(0)&=\left(\left(\frac br+
       b\left(\frac 1m-1\right)r-\frac{b^{2m+1}r^{1-2m}}m\right)e^{it},
         \sqrt{(1-b^{2m}r^{-2m})(1-b^{2m}r^{2-2m})}e^{is}\right)\\
     &=:(\gamma_1(r)e^{it},\gamma_2(r)e^{is}).
    \end{aligned}$$
The mapping
\begin{equation}\label{map}
\Delta\times[0,2\pi)\times(b,1)\ni(\zeta,t,r)\longmapsto
      \zeta(\gamma_1(r)e^{it},\gamma_2(r))
\end{equation}
parametrizes $I_{12}$. We will need a lemma.

\begin{lem} \label{lem}
Let $F(\zeta,z)=\zeta(f(z),g(z))$ be a function of two complex variables, 
where $f$ and $g$ are $C^1$. 
Then the real Jacobian of $F$ is equal to $|\zeta|^2H(z)$, where
  $$H=|f|^2(|g_{\bar z}|^2-|g_z|^2)+|g|^2(|f_{\bar z}|^2-|f_z|^2)
    +2\Re\big(f\bar g(\overline{f_z}g_z
         -\overline{f_{\bar z}}g_{\bar z})\big).$$
\end{lem}

The proof is left to the reader. For the mapping \eqref{map} we can
compute that
  $$\begin{aligned}H&=\gamma_1\gamma_2
               (\gamma_1\gamma_2'-\gamma_1'\gamma_2)\\
    &=-\frac{b^2}{m^2}r^{-6 m-3}
        \left[b^{2m}\left(-mr^2+m-1\right)+r^{2m}\right]
        \left[r^{2m}\left((m-1)r^2+m\right)-(2m-1)r^2b^{2 m}\right]\\
     &\ \ \ \ \ \ \ \ \times\left[r^2 b^{2m}+r^{2 m} 
         \left((m-1)r^2-m\right)\right]
   \end{aligned}.$$
\begin{comment}
  $$\begin{aligned}H&=\gamma_1\gamma_2(\gamma_1\gamma_2'-\gamma_1'\gamma_2)\\
    &=\frac{b^{6m+2}r^{1-6m}}{m^2}-\frac{3b^{4m+2}r^{1-4m}}{m^2}
       +\frac{3b^{2m+2}r^{1-2m}}{m^2}+2b^{6m+2}r^{1-6m}\\
    &\ \ \ -\frac{3b^{6m+2}r^{1-6m}}{m}-2b^{6m+2}r^{3-6 m}
        +\frac{b^{6m+2}r^{3-6m}}{m}+2b^{4m+2}r^{-4m-1}\\
    &\ \ \ -2mb^{4m+2}r^{-4m-1}-6b^{4m+2}r^{1-4m}+4mb^{4m+2}r^{1-4m}
        +\frac{8b^{4m+2}r^{1-4m}}{m}\\
    &\ \ \ +4b^{4m+2}r^{3-4m}-2mb^{4m+2}r^{3-4m}-\frac{2b^{4m+2}r^{3-4m}}{m}
        -b^{2m+2}r^{-2m-3}\\
    &\ \ \ +mb^{2m+2}r^{-2m-3}-2b^{2m+2}r^{-2m-1}-mb^{2m+2}r^{-2m-1}
         +5b^{2m+2}r^{1-2m}\\
    &\ \ \ -mb^{2m+2}r^{1-2m}-\frac{7b^{2m+2}r^{1-2m}}{m}
         -2b^{2m+2}r^{3-2m}+mb^{2m+2}r^{3-2m}\\
    &\ \ \ +\frac{b^{2m+2}r^{3-2m}}{m}-\frac{b^2r}{m^2}+\frac{2b^2r}{m}
          +\frac{b^2}{r^3}-b^2r.
   \end{aligned}.$$
\end{comment}
Since 
\begin{equation}\label{di2}
\int_\Delta|\zeta|^2d\lambda(\zeta)=\frac\pi 2,
\end{equation}
we obtain
\begin{equation}\label{i12}
\begin{aligned}\lambda(I_{12})&=\pi^2\int_b^1|H|\,dr\\
    &=\pi^2\left(\frac{(1-2m)^2}{m^2(3m-1)(3m-2)}b^{6m+2}
          -\frac{3}{m^2(m+1)(m-2)}b^{2m+2} 
          -\frac{3}{2m^2}b^{4 m+2}\right.\\
    &\ \ \ \ \ \ \ \ \ \ \ \ \left.+\frac{m}{2(m-2)(3m-2)}b^6
          +\frac{3m}{3m-1}b^4-\frac{4m^2-m+1}{2m^2}b^2
        +\frac{m}{m+1}\right).
    \end{aligned}
\end{equation}

To compute the volume of $I_2$ we consider geodesics
of the form
  $$\varphi(\zeta)=\left(a_1\left(\frac{1-\bar\alpha_1\zeta}
     {1-\bar\alpha_0\zeta}\right)^{1/m},
         a_2\frac{\zeta-\alpha_2}{1-\bar\alpha_0\zeta}\right),$$
where $a_1,a_2\in\mathbb C_\ast$, $\alpha_0,\alpha_2\in\Delta$, 
$\alpha_1\in\bar\Delta$ satisfy \eqref{c1}, \eqref{c2}.
By \eqref{w0} and since $\varphi(0)=(b,0)$ we have $a_1=b$, 
$\alpha_2=0$ and by \eqref{c1} $\alpha_0=b^{2m}\alpha_1$. 
By \eqref{c2}
  $$1+b^{4m}|\alpha_1|^2=b^{2m}\left(1+|\alpha_1|^2\right)+|a_2|^2,$$
that is
  $$|a_2|^2=(1-b^{2m})(1-b^{2m}|\alpha_1|^2).$$
This means that any $\alpha_1\in\Delta$ is allowed and by \eqref{w1}
  $$\begin{aligned}\varphi'(0)
     &=\left(\frac{b(b^{2m}-1)}m\bar\alpha_1,a_2\right)\\
     &=\left(\frac{b(1-b^{2m})r}me^{it},\sqrt{(1-b^{2m})(1-b^{2m}r^2)}
      e^{is}\right),
     \end{aligned}$$
where $\alpha_1=-re^{-it}$, $a_2=\rho e^{is}$. Similarly as before
we have
  $$H=-\frac{b^2(1-b^{2m})^3r}{m^2}$$
and
  $$\lambda(I_2)=\pi^2\int_0^1|H|\,dr=\frac{\pi^2b^2(1-b^{2m})^3}{2m^2}.$$
This combined with \eqref{i12} finishes the proof of Theorem \ref{thm1}. \qed

\section{\bf Proof of Theorem \ref{thm2}}

For $\Omega=\mathcal E(1/2,1/2)$ and $w=(b,b)$, where $0<b<1/2$, we have
by \eqref{w0}
\begin{equation}\label{ww0}
  a_j=\begin{cases}\dis -\frac b{\alpha_j},\ &j\in A\\
      \dis b, &j\notin A\end{cases}
\end{equation}
and by \eqref{w1}
\begin{equation}\label{ww1}
  \varphi_j'(0)
     =\begin{cases}\dis 2b\bar\alpha_0-b\left(\bar\alpha_j
           +\frac 1{\alpha_j}\right),\ 
            &j\in A\\ \\
      \dis 2b(\bar\alpha_0-\bar\alpha_j),
              &j\notin A\end{cases}.
\end{equation}
There are four possibilities
for the set $A$: $\emptyset$, $\{1\}$, $\{2\}$, and $\{1,2\}$. 
Denote the corresponding parts of $I_\Omega(w)$ by $I_0$, $I_1$,
$I_2$, and $I_{12}$, respectively, so that
\begin{equation}\label{sum}
\begin{aligned}
\lambda(I_\Omega(w))
     &=\lambda(I_0)+\lambda(I_1)+\lambda(I_2)+\lambda(I_{12})\\
     &=\lambda(I_0)+2\lambda(I_1)+\lambda(I_{12}).
\end{aligned}
\end{equation}

\vskip 10pt

\noindent {\bf The case $A=\{1,2\}$}

\vskip 7pt

By \eqref{c1}, \eqref{c2} and \eqref{ww0}
\begin{equation}\label{eq1}
\left(\frac 1b+2b\right)|\alpha_1|\,|\alpha_2|+2b\Re(\alpha_1
\bar\alpha_2)=(1+|\alpha_1|^2)|\alpha_2|+
(1+|\alpha_2|^2)|\alpha_1|.
\end{equation}
Since the set of $\alpha\in\Delta^2$ satisfying \eqref{eq1}
is $S^1$-invariant, let us consider only those $\alpha$ with 
$\alpha_2>0$. If we then replace $\alpha_1$ with $\bar\alpha_1$
then \eqref{eq1} will still be valid and $\varphi'(0)$ will
be replaced by $\overline{\varphi'(0)}$. We thus consider
\begin{equation}\label{polar}
\alpha_1=re^{it},\ \ \ \alpha_2=\rho,\ \ \ r,\rho\in(0,1),
\ t\in(0,\pi);
\end{equation}
to get $\lambda(I_{12})$ we will have to multiply the obtained 
volume by 2. The condition \eqref{eq1} transforms to
\begin{equation}\label{eqa12}
\frac 1b+2b(1+\cos t)=r+\frac 1r+\rho+\frac 1\rho.
\end{equation}
It will be convenient to substitute $x=r+1/r$, $y=t$,
and consider the domain
\begin{equation}\label{dom}
U:=\big\{(x,y)\in(2,\frac 1b+4b-2)\times(0,\pi):
       x<\frac 1b+2b(1+\cos y)-2\big\}.
\end{equation}
We have
  $$\alpha_0=b\left(\frac{\alpha_1}{|\alpha_1|}+
     \frac{\alpha_2}{|\alpha_2|}\right)
     =b(e^{it}+1)$$
and thus by \eqref{ww1} and \eqref{eqa12}
\begin{equation}\label{m1}
\begin{aligned}
\varphi'(0)&=b\left(2\bar\alpha_0-\bar\alpha_1-\frac 1{\alpha_1},
    2\bar\alpha_0-\bar\alpha_2-\frac 1{\alpha_2}\right)\\
  &=\left(2b^2(e^{-it}+1)-b(r+\frac 1r)e^{-it},
      2b^2(e^{-it}+1)-b(\rho+\frac 1\rho)\right)\\
  &=\big(2b^2+b(2b-x)e^{-iy},bx-1-2b^2i\sin y\big)\\
  &=:(f(z),g(z)).
\end{aligned}
\end{equation}
The mapping
  $$\Delta\times U\ni(\zeta,z)\longmapsto\zeta(f(z),g(z))$$
parametrizes $I_{12}$. From Lemma \ref{lem} and \eqref{di2}
it follows that 
  $$\lambda(I_{12})=\pi\iint_U|H|d\lambda,$$
where $f$, $g$ are given by \eqref{m1}, $U$ by \eqref{dom} 
(recall that again we had to multiply by 2) and we can compute that
  $$H=b^2\left[1-2 b^2(\cos y+1)\right]
    \left[-bx^2+(1+2 b^2(\cos y+1))(x-2b)
              -2b(b^2\cos(2y)+1)\right].$$
One can check that $H>0$ in $U$. The region $U$ may look as follows 

\vskip 10pt

\noindent\phantom{aa}
\includegraphics[scale=0.5]{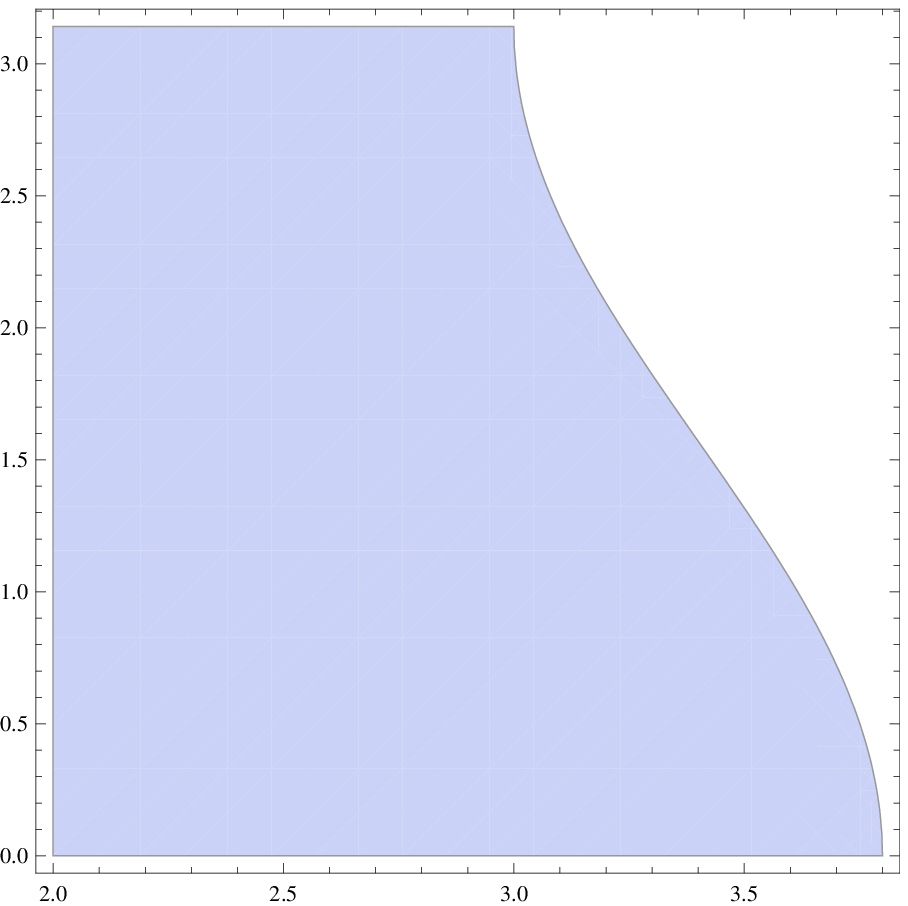}\hfill
\includegraphics[scale=0.5]{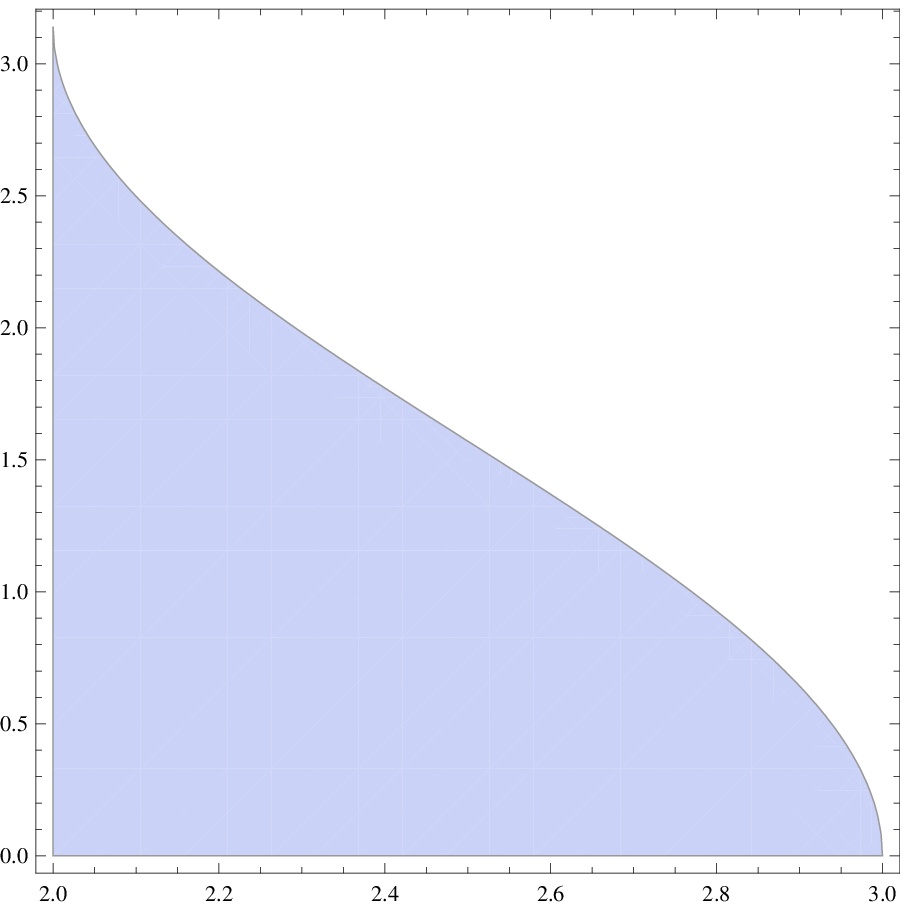}\hfill
\includegraphics[scale=0.5]{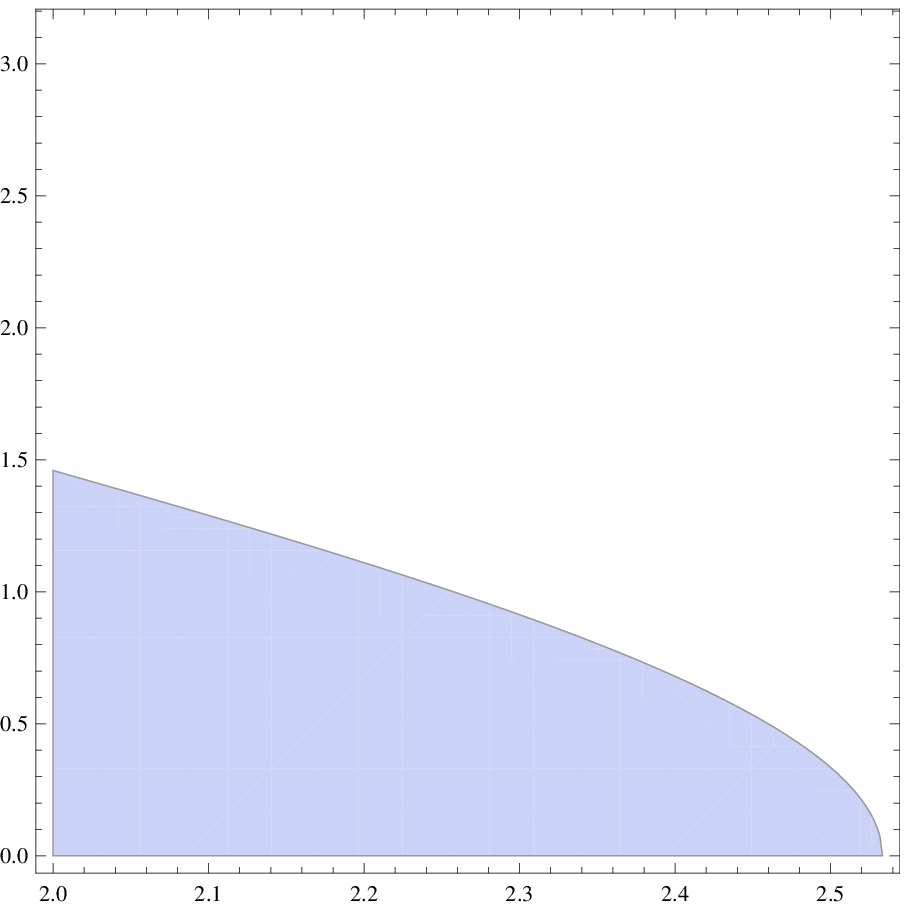}
\phantom{aa}

\vskip -10pt

\phantom{a}\hskip 20mm $b=0.2$ \hskip 40.5mm $b=0.25$
\hskip 40.5mm $b=0.3$

\vskip 7pt

\noindent We set 
  $$y_0:=\begin{cases}\pi&b\leq 1/4\\
    \dis\arccos\left(-1+\frac{4b-1}{2b^2}\right)\ &b>1/4
    \end{cases},$$
then
  $$\lambda(I_{12})=\pi\int_0^{y_0}\int_2^{1/b+2b(1+\cos y)-2}
     H\,dxdy.$$
For $b\leq 1/4$ we will get
\begin{equation}\label{i12-1}
\lambda(I_{12})=\frac{\pi^2}{6}
  (1 - 32 b^2 + 80 b^3 - 12 b^4 - 112 b^5 + 176 b^6 - 192 b^7 + 110 b^8)
\end{equation}
and for $b>1/4$
\begin{equation}\label{i12-2}
\begin{aligned}
  \lambda(&I_{12})=\frac{\pi}{72}(37-140b+270b^2-528b^3+530b^4
               -712b^5+660b^6)(1-2b)\sqrt{4b-1}\\
   &\ \ +\frac\pi 6(1-32b^2+80b^3-12b^4-112b^5+176b^6- 
      192b^7+110b^8)\arccos\left(-1+\frac{4b-1}{2b^2}\right).
\end{aligned}\end{equation}

\vskip 10pt

\noindent{\bf The case $A=\{1\}$}

\vskip 7pt

By \eqref{ww0} $a_1=-b/\alpha_1$, $a_2=b$ and by \eqref{c1}
$\alpha_0=b(\alpha_1/|\alpha_1|+\alpha_2)$. From \eqref{c2} 
we get 
\begin{equation}\label{al1}
  1+b^2\left(1+\frac{2\Re(\alpha_1\bar\alpha_2)}{|\alpha_1|}
       +|\alpha_2|^2\right)
    =\frac b{|\alpha_1|}(1+|\alpha_1|^2)+b(1+|\alpha_2|^2).
\end{equation}
We may assume that $\alpha_1>0$, then \eqref{al1} has 
a solution $\alpha_1\in(0,1)$ if and only if $T>2$, where
  $$\begin{aligned}
    T&=\frac 1b+b\big(1+2\Re\alpha_2+|\alpha_2|^2\big)
       -1-|\alpha_2|^2\\
     &=\frac 1b+b-1+2bx-(1-b)(x^2+y^2),
    \end{aligned}$$
and we write $\alpha_2=x+iy$. This means that
\begin{equation}\label{disk}
  \left|\alpha_2-\frac b{1-b}\right|<\frac{1-2b}{\sqrt b(1-b)}
\end{equation}
and the set $U$ will be the intersection of this disc with $\Delta$.
By \eqref{ww1} and \eqref{al1}
  $$\varphi'(0)=2b\big(b(1+\bar\alpha_2)-T/2,b-(1-b)\bar\alpha_2\big)$$
and therefore 
\begin{align*}
  f&=2b^2(1+x)-bT-2b^2yi,\\
  g&=2b^2-2b(1-b)x+2b(1-b)yi.
\end{align*}
We can compute that
  $$\begin{aligned}
  H&=4(1-b)b^2\big[b^2(1+2x)-(1-b)\big(1+b(x^2+y^2)\big)\big]\\ 
    &\ \ \ \ \ \ \ \ \ 
       \big[-1+2b+b^3-2b^2(1-b)x+b(1-b)^2(x^2+y^2)\big]\\
    &=4(1-b)b^3(b+b^2-(1-b)T)(b^2+2b-2+bT).
     \end{aligned}$$
One can check that $H>0$ everywhere on $U$.

If $b\leq 1/4$ then $U=\Delta$ and using the 
polar coordinates in $\Delta$ and Lemma \ref{lem} we will get 
\begin{equation}\label{i1}
 \lambda(I_1)=\frac{2\pi^2}3(1-b)b^2(3-9b+2b^2+6b^3-6b^4+10b^5).
\end{equation}
For $b>1/4$ it is more convenient to use the polar coordinates
in the disk \eqref{disk} instead:
\begin{equation*}\label{polar2}
  x=\frac b{1-b}+r\cos t,\ \ \ y=r\sin t,
\end{equation*}
then 
  $$H=4b^2(1-2b)^2-4b^4(1-b)^4r^4.$$
For $r$ with 
  $$\frac{1-2b}{1-b}<r<\frac{1-2b}{\sqrt b(1-b)}$$
the circles $\{|\alpha_2-b/(1-b)|=r\}$ and $\{|\alpha_2|=1\}$
intersect when $t=\pm t(r)$, where
\begin{equation}\label{tr}
  t(r)=\arccos\frac{1-2b-(1-b)^2r^2}{2br(1-b)}.
\end{equation}
Therefore
 $$\lambda(I_1)=\pi^2\int_0^{(1-2b)/(1-b)}rHdr
    +\pi\int_{(1-2b)/(1-b)}^{(1-2b)/(\sqrt b(1-b))}r(\pi-t(r))Hdr.$$
We can compute the second integral using the following 
indefinite integrals:
\begin{equation}\label{indef}
  \begin{aligned}
     \int v&\arccos\left(\frac av-v\right)dv
                   =\frac 14\sqrt{-a^2+2av^2-v^4+v^2}\\
     &\ \ +\frac{4a+1}8
        \arctan\frac{2a-2v^2+1}{2\sqrt{-a^2+2av^2-v^4+v^2}}
       +\frac{v^2}2\arccos\left(\frac av-v\right)+const,\\
     \int v^5&\arccos\left(\frac av-v\right)dv
          =\frac 1{288}\big(15+78a+80a^2+(10+32a)v^2+8v^4\big)
                  \sqrt{-a^2+2av^2-v^4+v^2}\\
     &\ \ +\frac{5+36a+72a^2+32a^3 }{192}
        \arctan\frac{2a-2v^2+1}{2\sqrt{-a^2+2av^2-v^4+v^2}}
       +\frac{v^6}6\arccos\left(\frac av-v\right)+const.
     \end{aligned}
\end{equation}
We will obtain
\begin{equation}\label{i1-2}
\begin{aligned}\lambda(I_1)
  &=-\frac{\pi ^2 b \left(10 b^9-36 b^8+54 b^7+84 b^6-375 b^5+414 b^4-166 b^3-6 b^2+21
     b-4\right)}{3(1-b)^2}\\
    &\ \ 
    +\frac{\pi b (1-2b)\left(30 b^6-58 b^5+43 b^4-19 b^3-26 b^2+32
       b-8\right)}{9 (1-b)}\sqrt{4b-1}\\
    &\ \
    +\frac{4\pi(1-2b)^4b\left(2 b^2-2 b-1\right)}{3(1-b)^2}
       \arccos\frac{3b-1}{2b^{3/2}}\\
    &\ \
    +\frac{2}{3}\pi(1-b)b^2\left(10 b^5-6 b^4+6 b^3+2 b^2-9 b+3\right) 
       \arctan\frac{2b^2-4b+1}{(1-2b)\sqrt{4 b-1}}.
\end{aligned}
\end{equation}
%\begin{equation}\label{i1-2}
%\begin{aligned}\lambda(I_1)
%  &=-\frac{2\pi b\left(10 b^9-36 b^8+54 b^7-108 b^6+201 b^5-162 b^4
%              +26 b^3+30 b^2-15b+2\right)}{3(1-b)^2}\\
%    &\ \ \ \ \ \ \ \ \ \ \ 
%         \times\arccos\left(-1+\frac{4 b-1}{2b^2}\right)\\
%    &\ \ 
%    +\frac{\pi b(1-2b)\left(30 b^6-58 b^5+43 b^4-19 b^3-26 b^2+32
%     b-8\right)}{9(1-b)}\sqrt{4b-1}\\
%    &\ \
%    +\frac{4\pi(1-2b)^4b\left(1+2b-2b^2\right)}{3(1-b)^2}
%       \left(\pi+\arctan\sqrt{4b-1}\right)
%\end{aligned}
%\end{equation}
for $b>1/4$.

\vskip 10pt

\noindent{\bf The case $A=\emptyset$}

\vskip 7pt

We have $a_1=a_2=b$ and $\alpha_0=b(\alpha_1+\alpha_2)$. Therefore
\begin{equation}\label{l01}
-b(1-b)(|\alpha_1|^2+|\alpha_2|^2)+2b^2\Re(\alpha_1\bar\alpha_2)
  +1-2b=0.
\end{equation}
Again, we may assume that $\alpha_1>0$. We may also assume that
$\Re\alpha_2\geq 0$ and then multiply the resulting integral by 2.
The equation \eqref{l01} has a solution $\alpha_1$ if  
  $$D:=-b(1-b)^2|\alpha_2|^2+b^3(\Re\alpha_2)^2+(1-b)(1-2b)
      \geq 0.$$
It satisfies $\alpha_1<1$ if
  $$Q:=\frac{b^{3/2}\Re\alpha_2+\sqrt{D}}{\sqrt b(1-b)}<1.$$
This means that
\begin{equation}\label{disk2}
  \left|\alpha_2-\frac b{1-b}\right|>\frac{1-2b}{\sqrt b(1-b)}.
\end{equation}
By $U$ we will denote the set of $\alpha_2\in\Delta$ satisfying
\eqref{disk2}. For $b\leq 1/4$ we have $U=\emptyset$ and thus 
$\lambda(I_0)=0$ then. This together with 
\eqref{sum}, \eqref{i12-1} and \eqref{i1} gives \eqref{b14}.

Assume that $b>1/4$. By \eqref{ww1}
  $$\varphi'(0)=2b\big((b-1)Q+b\bar\alpha_2,
          bQ+(b-1)\bar\alpha_2\big),$$
so that, 
\begin{align*}
  f&=2b\big((b-1)Q+bx\big)-2b^2y\,i\\
  g&=2b\big(bQ+(b-1)x)+2b(1-b)y\,i.
\end{align*}
One can compute that
  $$H=\frac{16b^3(1-2b)^3}{1-b}\left(1+\frac{b^{3/2}x}
    {\sqrt{D}}\right).$$
By Lemma \ref{lem}
  $$\lambda(I_0)=\pi\int_{-1}^{-1+(4b-1)/(2b^2)}
        \int_{y_2(x)}^{\sqrt{1-x^2}}Hdydx,$$
where
  $$y_2(x)=
    \begin{cases}
     0, &\ \ \ \ \ \ \ \ \ \ \ \ \ \ \!\!
      \dis -1\leq x\leq\frac{b^{3/2}+2b-1}{\sqrt b(1-b)}\\
      \sqrt{\dis\frac{(1-2b)^2}{b(1-b)^2}
        -\left(x-\frac b{1-b}\right)^2},\ \ \ \
      &\dis\frac{b^{3/2}+2b-1}{\sqrt b(1-b)}\leq x\leq 
          \dis-1+\frac{4b-1}{2b^2}.
   \end{cases}$$

It is clear from this formula that $\lambda(I_0)$ is analytic
for $b\in(1/4,1/2)$. We may therefore restrict ourselves
to the interval $(1/4,(3-\sqrt 5)/2)$, then $0\notin U$
and we will use polar coordinates in $\Delta$, that is
  $$x=r\cos t,\ \ \ \ y=r\sin t.$$
We will get
  $$\begin{aligned}
    \lambda(I_0)&=\frac{16\pi b^3(1-2b)^3}{1-b}\int_{r_0}^1
        r\left(\arccos\frac{1-3b+b^2-b(1-b)r^2}{2b^2r}\right.\\
      &\ \ \ \ \ \ \ \ \ \ \ \ \ \ \ \ \ \ \ \ \ \ \ \ \ \ \ 
         \ \ \ \ \ 
        \left.-\arctan\frac{\sqrt{4b^4r^2-(1-3b+b^2-b(1-b)r^2)^2}}
          {1-b-b^2-b(1-b)r^2}\right)dr,
     \end{aligned}$$
where 
  $$r_0=\frac{1-2b-b^{3/2}}{\sqrt b(1-b)}.$$
Using \eqref{indef} one can compute that 
  $$\begin{aligned}\int_{r_0}^1r&\arccos\frac{1-3b+b^2-b(1-b)r^2}{2b^2r}dr
     =\frac{\pi\left(2 b^3-8 b^2+6 b-1\right)}{4 (b-1)^2 b}\\
     &-\frac 12\arccos\left(-1+\frac{4b-1}{2 b^2}\right)
       +\frac{1-2b}{4 b(1-b)}\sqrt{4b-1}
       +\frac{(1-2b)^2}{2b(1-b)^2}\arctan\frac{1-3 b}{(1-b)\sqrt{4b-1}}.
    \end{aligned}$$
On the other hand, since
  $$\begin{aligned}
    \int\frac 1{v^2}\arctan&\sqrt{-av^2+v-1}\,dv=
       \frac 1{2v}\sqrt{-av^2+v-1}-\frac 1v\arctan\sqrt{-av^2+v-1}\\
      &-\frac a2\arctan\frac{2a\sqrt{-av^2+v-1}}{-av-2a+1}
      +\frac{2a-1}4\arctan\frac{(v-2)\sqrt{-av^2+v-1}}{2av^2-2v+2}
      +const,
  \end{aligned}$$
we will obtain
  $$\begin{aligned}
     \int_{r_0}^1r\arctan&\frac{\sqrt{4b^4r^2-(1-3b+b^2-b(1-b)r^2)^2}}
               {1-b-b^2-b(1-b)r^2}dx
            =\frac{\pi(1-2b)(b+1)}{8(1-b)^2}
        +\frac{1-2b}{4b(1-b)}\sqrt{4b-1}\\
     &-\frac{(b+2)(1-2b)}{4b(1-b)}\arctan\sqrt{4b-1}
    -\frac{(1+b)(1-2b)}{4(1-b)^2}\arctan\frac{1-3 b}{(1-b)\sqrt{4b-1}}.
  \end{aligned}$$
Therefore
\begin{equation}\label{i0-2}
  \begin{aligned}\lambda&(I_0)=
     \frac{2\pi^2b^2(1-2b)^3(-6b^2+9 b-2)}{(1-b)^2}
     -\frac{8\pi b^3(1-2b)^3}{1-b}
                \arccos\left(-1+\frac{4b-1}{2 b^2}\right)\\
    &+\frac{4\pi b^2(1-2b)^4(b+2)}{(1-b)^2}\arctan\sqrt{4b-1}
      +\frac{4\pi b^2(1-2 b)^4(2-b)}{(1-b)^2}
          \arctan\frac{1-3 b}{(1-b)\sqrt{4b-1}}.
 \end{aligned}\end{equation}
Using the formulas
\begin{equation}\label{arccos}
\arccos\left(-1+\frac{4b-1}{2b^2}\right)
     =\arctan\frac{2b^2-4b+1}{(1-2b)\sqrt{4b-1}}+\frac\pi 2,
\end{equation}
  $$\arccos\frac{3b-1}{2b^{3/2}}=\arctan\sqrt{4b-1}
       -\arctan\frac{2b^2-4b+1}{(1-2b)\sqrt{4b-1}}+\frac\pi 2,$$
and combining \eqref{sum}, \eqref{i12-2}, \eqref{i1-2} and \eqref{i0-2},
we get \eqref{a14} for $b>1/4$.

Denoting by $\chi_-$ and $\chi_+$ the functions defined by the 
right-hand sides of \eqref{b14} and \eqref{a14}, respectively,
we can compute that at 1/4
  $$\chi_-=\chi_+=\frac{15887}{196608}\pi^2,
  \ \ \ \chi'_-=\chi'_+=-\frac{3521}{6144}\pi^2,
  \ \ \ \chi''_-=\chi''_+=-\frac{215}{1536}\pi^2,
   \ \ \ \chi^{(3)}_-=\chi^{(3)}_+=\frac{1785}{64}\pi^2,$$
but
  $$\chi^{(4)}_-=\frac{1549}{16}\pi^2,\ \ \ \ \ \
    \chi^{(4)}_+=\infty.$$
This shows that our function is $C^3$ but not $C^{3,1}$ at 
$1/4$. \qed

In fact, using \eqref{arccos} and
  $$\arctan(1/x)=\frac\pi 2-\arctan x,\ \ \ x>0,$$
for $b\in (1/4,1-1/\sqrt 2)$ the formula \eqref{a14} can be written as
\begin{equation*}
\begin{aligned}
  \lambda(I_\Omega&((b,b)))=
    \frac{\pi^2}6
          \big(30b^8-64b^7+80b^6-80b^5+76b^4-16b^3-8b^2+1\big)\\
    &+\frac{\pi(1-2 b)\left(-180 b^7+444 b^6-554 b^5+754 b^4-1214 b^3+922
       b^2-305 b+37\right)}{72 (1-b)}\sqrt{4 b-1}\\
    &+\frac{4\pi b(1-2b)^4\left(7 b^2+2 b-2\right)}{3(1-b)^2}
             \arctan\sqrt{4b-1}\\
    &+\frac{\pi\left(30 b^{10}-124 b^9+238 b^8-176 b^7-
     260 b^6+424 b^5-76 b^4-144 b^3+89b^2-18 b+1\right)}{6(1-b)^2}\\
    & \ \ \ \ \ \ \ \ \ \times\arctan\frac{(1-2 b)\sqrt{4 b-1}}{2 b^2-4 b+1}\\
    &-\frac{4\pi b^2(1-2 b)^4(2-b)}{(1-b)^2}
        \arctan\frac{(1-b) \sqrt{4b-1}}{1-3 b}.
\end{aligned}
\end{equation*}

\end{document}